\documentclass[11pt]{article}

\usepackage{epic,eepic} 
\usepackage{psfrag}
\usepackage{graphicx}
\graphicspath{{./figures/}} 
\usepackage[a4paper]{geometry}
\usepackage{subfigure}
\usepackage{makeidx}

\usepackage{brem_book}

\usepackage{pst-node}
\psset{xunit=1cm,yunit=1cm}

\def\beq#1\eeq{\begin{equation}#1\end{equation}}

\makeindex

\title{Sampling cluster point processes: a review}

\author{Pierre Br\'emaud\,\thanks{\'Ecole Polytechnique F\'ed\'erale de Lausanne, Switzerland}}

\date{}
\begin{document}


\maketitle

\sk\noi {\bf Abstract:} The theme of this article is the sampling of cluster and iterated cluster point processes. It is partially a review, mainly of the Brix--Kendall exact sampling method for cluster point processes and its adaptation by M{\o}ller and Rasmussen to Hawkes branching point processes on the real line with light-tail fertility rate. A formal proof via Laplace transforms of the validity of the method in terms of general clusters that are not necessarily point processes fits this purpose and allows to include the exact sampling of Boolean models. The main novel aspect of this review is the extension of the above sampling methods to non-Poissonian germ point processes. 

\section{Introduction} 

Sampling a probability distribution $Q$ on a measurable space $(E,\E)$ consists, by definition, to generate a random element $X$ of distribution is $Q$. For doing this, it is assumed that one has at disposition any number of copies of ''easily generated'' random elements, such as for instance, {\sc iid} random variables uniformly distributed on a unit interval or Poisson processes. Implicit in the definition of sampling is the necessity that the generation of $X$ should necessitate only a finite (random) number of operations. Sampling is called exact sampling in order to insist on the difference with approximate sampling, for instance via Monte Carlo methods. 

\sk\noi In the domain of point processes the distributions to be sampled are those of a point process with an almost-surely finite number of points, in general (as will be the case in the following) the restriction of a point process to some domain, called the window, such that there is an almost-surely finite number of points in it. 

\sk\noi The main issue resides in the range of interaction of the distribution of the original point process. By this, we mean that the absence or presence of a point in the window may depend on the position of points outside of the window, points that will have to be, in principle, generated, and that are {\it a priori} in infinite number.

\sk\noi In their seminal work, Brix and Kendall have shown how to avoid this difficulty. The basic observation is that only a finite number of points outside the window are actually responsible for points inside the window. In the case they studied, the distribution of these ''active'' external point process is identifiable and a sample of it is easy to obtain. Therefore, the method consists in generating this external point process and sample for each of them its effect in the window. 
The method was applied to cluster point processes with a germ point process that is Poisson, such as the Cox cluster point process. In fact, the active external point process is obtained by thinning the original germ point process with a thinning probability that depends only on the location of the point independently of the rest of the cluster point process. Therefore, the thinned germ point process is also a Poisson point process with, under mild conditions, a finite number of points. 

\sk\noi Implementation of Brix--Kendall's requires knowledge of the thinning probability function. Such function is not always available in closed form. This situation occurs for instance when the typical cluster is a Hawkes process. M{\o}ller and Rasmussen have shown how deal with this situation in the case where the typical cluster is a standard Hawkes point process with a light-tail condition on the fertility rate. 

\sk\noi The purpose of the present article is to review the above methods and propose various extensions. A proof of the validity of the Brix--Kendall's algorithm in terms of Laplace transforms, thus making the intuitive arguments rigorous and at the same time emphasizing the universal relevance of the method, in particular to the exact sampling of Boolean models on a bounded window, an important application since the statistics of the Boolean model are usually not computable. 

\sk\noi We then show how to apply the method to a germ point that is not Poisson. Three cases are treated. 

\sk\noi (a) The germ point process lies on a grid (say, $\setZ^2$ or part of it, or a deterministic germ point process, 

\sk\noi (b) In the univariate case, a renewal process, and more generally any easily sampled point process with a bounded stochastic intensity, and  

\sk\noi (c) The point process is 

\sk\noi The M{\o}ller--Rasmussen method is extended to a Poisson germ point process of the class (b).

\section{Cluster point processes} 

Recall a few definitions and notations. Let $E$ be a locally compact topological space with a denumerable 
base (for short, \emph{\color{blue} l.c.d.b.}). Let $\B(E)$ be the Borel sigma-field on this topological space, 
that is, the sigma-field generated by the open sets of the topology. A subset of $E$ is called relatively compact\index{compact!relatively ---} if its closure is compact.

\sk\noi Let $M(E)$ be the set of \emph{locally finite} measures (that is taking finite values on locally compact sets) on $(E,\B(E))$ and let $\M(E)$ be the sigma--field 
on $M(E)$ generated by the mappings $p_C:\,\mu\to\mu(C)$, $C\in\B(E)$. A measure $\mu\in M(E)$ taking integer (possibly infinite) values is called a \emph{\color{blue} point measure}\index{point measure}. Such a point measure can be represented as a countable sum of Dirac measures
$$
\mu =\sum_{n\in S(\mu)} \varepsilon_{x_n} \, ,
$$
where $S(\mu)$ is a subset of $\setN$ and the $x_n$'s need not be distinct. The subset of $M(E)$ consisting of the locally finite point measures is denoted by $M_p(E)$\index{MEc@$M_p(E)$}, and we define the sigma-field $\M_p(E)$\index{MEd@$\M_p(E) $} on it as the sigma-field generated by the collection of sets 
$\{\mu \in M_p(E); \, \mu(C)\in F\},\,  C\in \B(E),\, F\in \B(\overline{\setR}_+)$.

\sk\noi A \emph{\color{blue} locally finite point process} on $E$ is a measurable mapping $N:(\Omega, \F)\to (M_p(E), \M_p(E))$. It is called \emph{\color{blue} simple} if $P(N(\{x\})\leq 1\mbox{ for all }x\in E)=0$. The \emph{\color{blue} intensity measure} of $N$ is the measure $\nu$ defined by $\nu(C)=E\cro{N(C)}$ ($C\in\B(E)$)

\sk\noi Let $N_0$ be a simple locally finite point process on the l.c.d.b. space $E$, with sequence of points $\{X_{0,n}\}_{n\in\setN}$ and locally finite intensity measure $\nu_0$. Let $\{Z_n\}_{n\in\setZ}$ be an {\sc iid} sequence of random measurable kernels from $(E \times \Omega, \B(E)\otimes \F)$ to $(E, \B(E))$, independent of $N_0$, and such that $E\cro{Z_1(x, C}=K(x,C)$ for a measurable kernel $K$ from $(E , \B(E))$ to $(E, \B(E))$ such that for all bounded $C\in \B(E)$
$$
\int_E K(x, C-x) \nu_0 (dx)<\infty \, .
$$
The random measure $N$ on $E$ defined by  
\beq \label{cluster1}
N(C)\defeq\sum_{n\in\setN} Z_n(X_{0,n} ,C-X_{0,n})
\eeq
is called a \emph{\color{blue} cluster random measure} with \emph{\color{blue} germ} $N_0$. The random measure $Z_n(X_{0,n} ,\cdot -X_{0,n})$ is the \emph{\color{blue} cluster} at $X_{0,n}$.

\sk\noi A straightforward application of Campbell's formula gives for the intensity measure $\nu$ of $N$
\beq 
\label{cluster2}
\nu(C)=\int_{E} K(y, C-y)\, \nu_0(dy)  \, .
\eeq
In fact, 
\begin{align*}
E\cro{\sum_{n\in\setN} Z_n(X_{0,n}\, ,\,C-X_{0,n})}
&=E\cro{E\cro{\sum_{n\in\setN} Z_n(X_{0,n}\, ,\,C-X_{0,n}) \, | \, \F^{N_0}}}\\
&=E\cro{\sum_{n\in\setN} E\cro{Z_n(X_{0,n}\, ,\,C-X_{0,n}) \, | \, \F^{N_0}}}\\
&=E\cro{\sum_{n\in\setN} \nu_Z(X_{0,n}\, ,\,C-X_{0,n}) }=\int_{E}\nu_Z(y,C-y))\, \nu_0(dy) \, .
\end{align*}
If this measure is locally finite, the random measure considered is a random elements of $M(E)$.

\sk\noi  
When for all $x\in E$ and all $\omega$, $Z_1(x,\omega, \cdot)\in M_p(E)$, (\ref{cluster1}) defines a cluster point process. This point process is simple if, for instance, $\nu_0$ is a diffuse measure.

\sk\noi When $Z_1$ ''does not depend on $x$'', that is when it is a random measurable kernel from $(\Omega,  \F)$ to $(E, \B(E))$, we use the notation $N=N_0*Z$, where $Z$ stands for the ``generic'' cluster, that is any random measure with the common distribution of the $Z_n$'s. Implicit in this notation is the assumption that the marks $Z_n$ of $N_0$ are {\sc iid} and independent of $N_0$. In this case, 
\beq 
\label{cluster3}
\nu=\nu_0*\nu_{Z_1} \, , 
\eeq 
where $\nu_{Z_1}$ is the intensity measure of $Z_1$. 

\sk\noi Note that in the case of point process clusters, the $Z_n$'s may have a point at $0$ in which case some, or all, points of the germ point process are part of the cluster point process. When $E=\setR^m$, a sufficient condition for the cluster point process to be simple is that its intensity measure be diffuse. This is the case whenever one of the measures 
of the convolution (\ref{cluster3}) is a multiple of the Lebesgue measure, and the other is a finite measure. For instance, if the intensity measure of the germ point process is of the form $\nu_0(dx)=\lambda_0 \ell^m(dx)$, then  
$$
\nu(C)=\int_{\setR^m} \lambda_0 \ell^m(C-x)\nu_Z(dx) =\int_{\setR^m} \lambda_0 \ell^m(C)\nu_Z(dx) =\lambda_0 \nu_Z(E) \ell^m(C)\, .
$$
We leave to the reader the task of finding general conditions 
that make of $N$ a \emph{simple} point process.

\sk\noi If $Z_1(x, \cdot)$ is a Poisson process, the cluster point process is called a Cox cluster point process.

\subsection{Branching point process} 
The point process $Z$ is a called a branching point process with single ancestor point\index{branching point process!--- with a single ancestor} at $0$ if
$$
Z=z_0 +z_1 + z_2 + \cdots 
$$
where $z_0\defeq \varepsilon_0$ (the point process with a single point, at $0$), and for all $n\geq 0$, $z_{n+1}$ is the cluster point process with germ point process $z_n$ and typical cluster $\alpha$, a simple finite point process such that $\alpha (\{0\})=0$. We may use the notation introduced a few lines above: $z_{n+1}=z_n*\alpha$, but remember that in this notation, the underlying clusters of the $n$-th generation that are attached to the $n$-th generation germ point process $z_n$ are {\sc iid} and independent of $z_n$. Moreover the collection of clusters of all generations are {\sc iid}.

\sk\noi In particular, the sequence $\{\sum_{k=0}^nz_k(E)\}_{n\geq 0}$ is a Galton-Watson process with 
a single ancestor and typical progeny distributed as $\alpha (E)$. In particular, if 
$E\cro{\alpha (E)}<1$, $Z$ is a finite point process and $E\cro{Z(E) }=\frac{1}{1-E\cro{\alpha (E)}}<\infty$, or with the notation $|\alpha|\defeq \alpha(E)$, 
$$
E\cro{Z(E) }=\frac{1}{1-E\cro{|\alpha |}}<\infty \, .
$$

\sk\noi We now define a general branching process as a particular kind of cluster point process with the following specificities: 

\sk (a) $E=\setR^m$, 

\sk (b) the intensity measure of the germ point process 
is $\nu_0$, and 

\sk (c) $Z\defeq Z_1$ is the branching point process with ancestor point at $0$ just described, where it is assumed that $E\cro{\alpha (\setR^m)}<1$. The generic cluster $\alpha$ is called the generic \emph{\color{blue} progeny} of the branching cluster point process.

\sk\noi This point process is also called an iterated cluster process\index{cluster point process!iterated ---}\index{branching cluster point process}\index{point process!branching cluster ---}, since it consists in a succession of generations, $N_0, N_1, N_2, \ldots$, where for $n\geq 1$, 
$N_n$ is obtained by $\alpha$-clustering of $N_{n-1}$, that is $N_n=N_{n-1}*\alpha $. The final point process being 
$$
N=\sum_{n\geq 0} N_n \, .
$$
The intensity measures of the successive generations are 
$\nu_0, \nu_1, \nu_2, \ldots$, where $\nu_n=\nu_{n-1}*\nu_\alpha$. In particular, for all $n\geq 0$,
$$
\nu_n(C)=\int_{\setR^m}\nu_{n-1}(C-x) \, \nu_\alpha(dx) \, .
$$
In the special case where $\nu_0=\lambda_0 \ell^m$, $\nu_n=\lambda_n \ell^m$ where $\lambda_n=\lambda_0 |\nu_\alpha|^n$. 
Finally, the intensity measure of $N$ is, since $|\nu_\alpha| <1$
$$
\nu (dx)= \frac{\lambda_0}{1-|\nu_\alpha|}\ell^m  (dx)\, .
$$

\sk\noi \index{sample!approximate ---!--- of a point process distribution}There are cases where exact sampling is not possible. We then must have recourse to approximate samples. An approximate sample of a distribution $P_N$ is an exact sample of a  distribution $P_{\widetilde{N}}$ ``close'' to $P_N$, where the closeness is measured in terms of the variation distance 
$$
d_V (P_N, P_{\widetilde{N}})=sup_{\Gamma\in\M_p(E)}|P_N(\Gamma)- P_{\widetilde{N}}(\Gamma)|.
$$

\subsection{Approximate samples of a cluster point process}
Consider the problem of generating a sample of the branching process on a ``window'' $W\in \B(\setR_m)$ of finite Lebesgue measure. This requires to construct the branching processes attached to all the points of the germ process, which are possibly in infinite number. This is in general out of reach (cases where this is possible will be considered later on). For the time being, suppose that instead of $N$, one succeeds in sampling its approximation 
$$
N^{(n)}=\sum_{k= 0}^n N_k \, ,
$$
whose intensity measure is, assuming that the intensity measure is $\lambda_0 \ell^m$,  $\lambda_0 \frac{1-\nu_\alpha (\setR^m)^n}{1-|\nu_\alpha |}\ell^m  (dx)$. In particular 
$$
E\cro{(N-N^{(n)})(W) }=|\nu_\alpha |^n \frac{\lambda_0}{1-|\nu_\alpha |}\ell^m  (W) \, ,
$$
and therefore, since for any integer-valued random variable $Y$, $P(Y>0)\leq E\cro{Y}$, 
$$
P((N-N^{(n)})(W)>0 )\leq \gamma |\nu_\alpha |\, ,
$$
where 
$$
\gamma = \frac{\lambda_0}{1-|\nu_\alpha |}\ell^m  (W) \, .
$$
This says that the probability that $N^{(n)}\not\equiv N$ on $W$ is lesser that $\gamma \nu_\alpha (\setR^m)^n$. In still other words, denoting by $N_W$ the restriction of $N$ to $W$,
$$
d_V(N_W^{(n)},N_W)\leq \gamma |\nu_\alpha |^n \, ,
$$
where $d_V$ is the variation distance. Indeed:
\begin{align*}
&|P(N_W^{(n)}\in \Gamma)- P(N_W\in \Gamma)|\\
=&|P(N_W^{(n)}\in \Gamma, N_W^{(n)}\equiv N_W)+P(N_W^{(n)}\in \Gamma, N_W^{(n)}\not\equiv N_W) \\
&\hphantom{aaaaaaa} - P(N_W\in \Gamma, N_W^{(n)}\equiv N_W)-P(N_W\in \Gamma, N_W^{(n)}\not\equiv N_W)| \\
&=|P(N_W\in \Gamma, N_W^{(n)}\equiv N_W)+P(N_W^{(n)}\in \Gamma, N_W^{(n)}\not\equiv N_W)\\
&\hphantom{aaaaaaa} - P(N_W\in \Gamma, N_W^{(n)}\equiv N_W)-P(N_W\in \Gamma, N_W^{(n)}\not\equiv N_W)| \\
&=|P(N_W^{(n)}\in \Gamma, N_W^{(n)}\not\equiv N_W)-P(N_W\in \Gamma, N_W^{(n)}\not\equiv N_W)|\leq P(N_W^{(n)}\not\equiv N_W) \, .
\end{align*}

\sk\noi Still, the above mentioned difficulty remains. However, if the support of $\alpha$ is finite, that is, if for some $R<\infty$, $P(\alpha (\{x\in\setR^m \, ; \, ||x||\geq R\})=0)=1$, a little thought shows that to obtain $N^{(n)}$ on $W$, it suffices to construct the branching processes attached to only the germ points at a distance less that $nR$ from $W$. These are in finite number. We therefore obtain an approximation of the sample we looked for, but the quality of this approximation in terms of the variation distance can be controlled, and made as good as desired by a proper choice of $n$.

\section{The Brix--Kendall's result via Laplace transforms} 

Consider the problem of sampling the distribution of the restriction of a cluster point process $N$ to $W\subset E$, denoted by $N_W$. As we noted before, this requires in principle to generate all the points $X_{0,n}$ of $N_0$, since all the associated point processes $Z_n(\cdot -X_{0,n})$ are suceptible to produce points in $W$. This is not feasible if there is an infinite number of points of the germ point process $N_0$ outside the window $W$. One solution is to approximate $N_W$ by taking into account only the points of the germ point process that are in a ``sufficiently'' large window $W'\supset W$. This introduces egde effects, here a loss of points in the window $W$. This ailment found a remedy in the case of the Cox cluster point process with the Brix--Kendall exact sampling algorithm whose natural idea is the following. All points of the germ point process do not contribute to $N_W$. Such a point located at $x$ will contribute only its cluster has at least a point in $W$, which happens with probability $1-e^{K(x,W-x)}$. The contributing germ points form a Poisson process of intensity measure $(1-e^{K(x,W-x)}) \, \widetilde{\mu}(dx)$, where $\widetilde{\mu}$ is the intensity measure of the germ point process. Therefore, it is reasonable to obtain a sample of $N_W$ by replacing the germ point process by a Poisson process of intensity measure $(1-e^{K(x,W-x)}) \, \widetilde{\mu}(dx)$, a finite point process if $\int_{\setR^d}(1-e^{x,K(W-x)}) \, \widetilde{\mu}(dx)<\infty$ which we henceforth suppose. For each point of this new germ point process, generate a cluster $\widetilde{N}_n$ whose distribution is that of any of the original clusters, only \emph{conditioned} by the event that it has at least one point in $W$. The sample of $N_W$ is then the sum of these new cluster point processes restricted to $W$.

\sk\noi The fundamental idea of the Brix--Kendall exact sampling algorithm is rather intuitive. It nevertheless requires a formal proof. This proof will be given in terms of Laplace transforms for random measures rather than just point processes. This generality will be used later for the exact sampling of Boolean models.

\sk\noi The computation of the Laplace transform of $N_W$ will prepare the way. Using the fact that $\{Z_n(X_{0,n}, \cdot)\}_{n\in \setZ}$ is, conditionally on $\F^{N_0}$, an independent sequence, 
\begin{align*}
L_{N_W}(\varphi)
&\defeq E\cro{\exp \acc{-\sum_{n\in\setN} \pth{\int_{W} \varphi(x) \, 
Z_n(X_{0,n}\,, \,dx-X_{0,n})}}}\\
&=E\cro{ E\cro{\exp \acc{-\sum_{n\in\setN} \pth{\int_{W} \varphi(x) \, 
Z_n(X_{0,n}\,, \,dx-X_{0,n})}}\, | \, \F^{N_0}}}\\
&=E\cro{E\cro{ \prod_{n\in\setZ}\exp \acc{-\int_{W} \varphi (x)\, 
Z_n(X_{0,n}\,, \,dx-  X_{0,n})   }\, | \, \F^{N_0}}    }\\
&=E\cro{\prod_{n\in\setN}E\cro{ \exp \acc{-\int_{W} \varphi (x)\, 
Z_n(X_{0,n}\,, \,dx-  X_{0,n})   }\, | \, \F^{N_0}}    }.
\end{align*} 
Now, with $A_n\defeq \{Z_n(x,W-X_{0,n})>0\}$,
\begin{align*}
&E\cro{ \exp \acc{-\int_{W} \varphi (x)\, Z_n(X_{0,n}\,, \,dx-  X_{0,n})   }\, | \, \F^{N_0}} \\
&\hphantom{aaaa}=E\cro{ \exp \acc{-\int_{W} \varphi (x)\, Z_n(X_{0,n}\,, \,dx-  X_{0,n})   }\, | \, X_{0,n}}\\
&\hphantom{aaaa}=E\cro{ \exp \acc{-\int_{W} \varphi (x)\, Z_n(X_{0,n}\,, \,dx-  X_{0,n})   } 1_{A_n}\, + 1-1_{A_n}| \, X_{0,n}} \\
&\hphantom{aaaa}=E\cro{ \exp \acc{-\int_{W} \varphi (x)\, Z_n(X_{0,n}\,, \,dx-  X_{0,n})   } 1_{\{Z_n(W-X_{0,n})>0\}}\, | \, X_{0,n}} \\
&\hphantom{aaaaaaaaaaaaaaaaaaaaaaaaaaaaaaaaaaaa}+ P(Z_n(X_{0,n}\,, \,W-X_{0,n})=0 \, | \, X_{0,n})\\
&\hphantom{aaaa}\defeq g_1(X_{0,n})+ g_2(X_{0,n}) \, .
\end{align*}
Therefore 
\begin{align*}
L_{N_W}(\varphi)&=E\cro{\prod_{n\in \setN}(g_1(X_{0,n})+ g_2(X_{0,n}))}\\
&=E\cro{e^{\sum_{n\in \setN}\log (g_1(X_{0,n})+ g_2(X_{0,n})) }   }\\
&=E\cro{e^{\int_E\log (g_1(x)+ g_2(x))  \, N_0(dx) }  }\\
&=E\cro{e^{\int_E\log (g_1(x)+ P(Z_1(x,W-x)>0))  \, N_0(dx)}   }\, , 
\end{align*} 
where $N_x$ is a typical cluster with germ at $x\in E$, that is, with the same distribution any of the $Z_n(\cdot -x)$. Now 
$$
g_1(x) +P(Z_1(x,W-x)>0)=\frac{g_1(x)}{P(Z_1(x,W-x)>0)}P(Z_1(x,W-x)>0) +P(Z_1(x,W-x)>0) \, .
$$
Observe that 
$$
\frac{g_1(x)}{P(Z_1(x,W-x)>0)}=E\cro{e^{-\int_W \varphi(y) \, Z_1(x,dy-x) }1_{\{Z_1(x,W-x)>0\}}    }/P(Z_1(x,W-x)>0)
$$
is the Laplace transform of the point process $Z_1(x,\cdot-x)$ conditioned to have at least one point in $W$.

\sk\noi We now turn to the Brix--Kendall exact sampling proposition, and construct a point process $\widetilde{N}_W$ on $W$ as follows. First, the point process $N_0$ is thinned, a point $X_{0,n}$ being retained with probability $p(X_{0,n})$ where $p(x)\defeq P(Z_1(x,W-x)>0)$ (defined above). More precisely, the thinned point process 
$\widetilde{N_0}$ is defined by
$$
\widetilde{N_0}(C)\defeq \sum_{n\in\setZ} 1_C(X_{0,n}) Y_n \, ,
$$ 
where $\{Y_n\}_{n\in\setZ}$ is, conditionally on $\F_{N_0}$, an independent sequence with values in $\{0,1\}$, and for each $n\in\setZ$,  $P(Y_n=1 \, | \, \F^{N_0})=P(Y_n=1 \, | \, X_{0,n})=p(X_{0,n})$. 
Then for each $n\in\setN$, let $\widetilde{Z}_n$ be a point process that has the same distribution as $Z_n$ conditioned by $Z_n(X_{0,n}\, , \,W -X_{0,n})>0$. This point process $\widetilde{Z}_n$ is obtained by sampling independent point processes of the type $Z_n$ until the condition $Z_n(X_{0,n}\, , \,W -X_{0,n})>0$ is satisfied. The candidate sample $\widetilde{N}$ is then constructed as 
$$
\widetilde{N}_W (C)= \sum_{n\in\setZ} Y_n\widetilde{Z}_n (X_{0,n}\, , \,C-X_{0,n}) \, .
$$
(Therefore, only the $\widetilde{Z}_n$ corresponding to a point $X_{0,n}$ that has been retained will need to be sampled.) In order to check that $\widetilde{N}_W$ is the desired exact sample, it must be proved that it has the same distribution as $N_W$. This is done below by showing that they have the same Laplace functional. 
Write 
\begin{align*}
L_{\widetilde{N}_W}(\varphi)
&\defeq E\cro{\exp \acc{-\sum_{n\in\setN} \pth{\int_{W}\varphi(x) \,Y_n\widetilde{Z}_n(X_{0,n}\, , \,dx-X_{0,n})}}}\\
&=E\cro{ \prod_{n\in\setZ}\exp \pth{-\int_{W} \varphi (x)\, Y_n\widetilde{Z}_n(X_{0,n}\, , \,dx-  X_{0,n})   }    }\\
&=E\cro{E\cro{ \prod_{n\in\setZ}\exp \acc{-\int_{W} \varphi (x)\, Y_n\widetilde{Z}_n(X_{0,n}\, , \,dx-  X_{0,n})   }\, | \, \F^{N_0}}    }\\
&=E\cro{\prod_{n\in\setN}E\cro{ \exp \acc{-\int_{W} \varphi (x)\, Y_n\widetilde{Z}_n(X_{0,n}\, , \,dx-  X_{0,n})   }\, | \, \F^{N_0}}    }\, .
\end{align*} 
Now write
\begin{align*}
&E\cro{
 \exp \acc{-\int_{W} \varphi (x)\, Y_n\widetilde{Z}_n(X_{0,n}\, , \,dx-  X_{0,n})   }\, | \, \F^{N_0}
} \\
&\hphantom{aaaa}
=E\cro{ \exp \acc{-\int_{W} \varphi (x)\, \widetilde{Z}_n(X_{0,n}\, , \,dx-  X_{0,n})}Y_n\, | \, X_{0,n}
} + 
E\cro{  1-Y_n\, | \, X_{0,n}} \\
&\hphantom{aaaa}=E\cro{ \exp \acc{-\int_{W} \varphi (x)\, \widetilde{Z}_n(X_{0,n}\, , \,dx-  X_{0,n})}\, | \, X_{0,n}} E\cro{  Y_n\, | \, X_{0,n}}+ 
E\cro{  1-Y_n\, | \, X_{0,n}} \\ 
&\hphantom{aaaa}=E\cro{ \exp \acc{-\int_{W} \varphi (x)\, \widetilde{Z}_n(X_{0,n}\, , \,dx-  X_{0,n})}\, | \, X_{0,n}}  p(X_{0,n}))+ (1-p(X_{0,n}))) \\ 
&\hphantom{aaaa}=g(X_{0,n})p(X_{0,n})+ 1-p(X_{0,n})\, , 
\end{align*}
where $g(x)$ is the Laplace functional of $Z_1(x, \cdot -x)$ conditioned to have at least one point in $W$. The rest of the verification is completed by 
\begin{align*}
L_{\widetilde{N}_W}(\varphi)
& E\cro{\prod_{n\in\setZ} (g(X_{0,n})  p(X_{0,n}) +1- p(X_{0,n})   }\\
&=E\cro{\exp\acc{\int_E \log (g(x)p(x)+1-p(x)) \, N_0(dx)     }     }\, ,
\end{align*}
and the observation $g(x)p(x)=g_1(x)$.

\sk\noi In the case where $N_0$ is a Poisson process of mean measure $\widetilde{\mu}$, the exact sampling procedure consits of constructing a thinned version $\widetilde{N}_0$ of $N_0$, in this case a Poisson process of intensity measure $P(Z_1(x,W-x)>0)\widetilde{\mu}(dx)$, and from each point $\widetilde{X}_{0,n}$ 
of $\widetilde{N}_0$ realize a version of $\widetilde{Z}_n$. There are two conditions for this to produce an exact sampling of $N_W$ in a finite number of operations. The number of points of $\widetilde{N}_0$ must be finite, a sufficient condition for this being that 
$$
\int_E P(Z_1(x,W-x)>0)\widetilde{\mu}(dx)<\infty \, .\eqno{(\star)}
$$
For of a Cox cluster point process $P(Z_1(x,W-x)>0)=1-e^{-K(x,W-x)}$ and 
therefore $\int_{\setR^d}P(Z_1(x,W-x)>0)\widetilde{\mu}(dx)=\int_{\setR^d}(1-e^{-K(x,W-x)})\, \widetilde{\mu}(dx) 
<\infty$ in view of condition ($\dagger$) and of the inequality $1-e^{-x} \leq x$ ($x\in\setR$).

\section{Exact sampling of Boolean models}

A Boolean set constructed on the germ point process $N_0$ is a random set of the form
$$
B\defeq \cup_{n\geq 1} (S_n +X_{0,n}) \, ,
$$
where the the sequence $\{S_n\}_{n\geq n}$ is an {\sc iid} sequence of closed random sets, for instance, closed balls centered at $0$ of radiuses $\{R_n\}_{n\geq n}$ forming an {\sc iid} sequence. The restriction of the Boolean set to the window $W$ is, by definition, the random set $B\cap W$. We make the assumption that almost-surely, $S_1$ is identical to the closure of its interior. Then, defining the random measure $Z_n$ by  
$$
Z_n(C)\defeq \int_C 1_{S_n}(x)\, dx =\ell^m (S_n\cap C)
$$
an exact sample of $B\cap W$ is obtained as soon as we have obtained an exact sample on the window $W$ of the random measure $N$ defined by
$$
N\defeq \sum_{n\geq 1} Z_n(\cdot -X_{0,n}) \, .
$$
The theory was done in sufficiently general terms to accomodate this case and to obtain the equally intuitive result that one should first thin the germ process with the thinning probability function
$$
p(x)\defeq P(\ell^m ((S_1 +x)\cap W)>0)\, .
$$
This is equivalent to 
\beq 
\label{boole}
p(x)\defeq P((S_1 +x)\cap W\neq \varnothing)\, .
\eeq
An important case that does not quite fit the above framework is when $S_1$ is a line passing through the origin $0$ (it is not identical to the closure of its interior). However, replacing the $S_n$'s by their fattened versions $S_n(\varepsilon)\defeq \{y\in\setR^m \, ;\, d(y,S_n)\leq \varepsilon\}$ fits the framework, and a limiting argument as $\varepsilon\to 0$ shows that the Brix--Kendall method applies with a thinning probability given by (\ref{boole}). 

\sk\noi For instance, suppose we seek to sample Poissonian lines inside a disk centered at $0$ and of radius $R$. By Poissonian lines we mean lines passing through the points of a homogeneous Poisson process on $\setR^2$ , say of intensity $\lambda$, with independent random uniform orientation. The probability of retaining a point of the germ point process located at $x$ is then 
$p(x)=\frac{1}{\pi}\arcsin \pth{\frac{R}{||x||}}$.

\section{Non-Poissonian germ processes}

Suppose that we take for granted that a typical cluster is easily generated. Then, as previous calculations confirm, there are two ingredients that make things work. First of all, the thinning probability function must be available in closed form. This is not the case for a Hawkes point process. This prompted M{\o}ller and Rasmussen to modify the thinning operation, as will be seen in the next section. 
A second possible difficulty when attempting to extend Brix--Kendall's method to non-Poissonian germ point processes is that to obtain a sample of the thinned germ point process. (In the case of a Poisson germ process, the thinned process is also a Poisson process and therefore this difficulty does not exist.)  

\sk\noi We now give three examples where the second limitation can be overcome. The first example, corresponding to extension (a) is an adaptation of the Poisson process case. Indeed, it consists in generating the random variable counting the number of points of the thinned point process, and then, place these points.

\sk\noi 
\subsection{Thinning the grid} Consider a point process on $\setN$ represented by a sequence $\{X_n\}_{n\geq 0}$ of {\sc iid} $\{0, 1\}$-valued random variables, with the common distribution given by 
$P(X_n=1)=p_n$ ($n\geq 0$). (We are therefore ``thinning the grid'' $\setN$, considered as a deterministic point process, with the thinning probability function $p_n$.) Suppose that $\sum_{n\geq 0}p_n<\infty$, which guarantees that the thinned grid has almost surely a finite number of points and let $T$ be its last point. Note that
$$
P(T=n)=P(X_n=1, X_{n+1}=0, X_{n+2}=0, \ldots)=p_n \prod_{k\geq n+1}(1-p_k) \eqno{(\star)}
$$
and that, for $0\leq k\leq n-1$,  
\begin{align*}
P(X_k=1 \mid T=n)&=\frac{P(X_k=1 ,T=n)}{P(T=n)}\\
&=\frac{P(X_k=1,X_n=1, X_{n+1}=0, X_{n+2}=0, \ldots)}{P(X_n=1, X_{n+1}=0, X_{n+2}=0, \ldots)}\\ 
&=\frac{P(X_k=1)P(X_n=1, X_{n+1}=0, X_{n+2}=0, \ldots)}{P(X_n=1, X_{n+1}=0, X_{n+2}=0, \ldots)}=P(X_k=1) \, .
\end{align*}
Therefore, in order to simulate the thinned grid, one may start by sampling a variable $T$ with the 
distribution ($\star$), and if $T=n$, set $X_n=1, X_{n+1}=0, X_{n+2}=0, \ldots$ and for  $0\leq k\leq n-1$, 
sample $X_k$ with the distribution $P(X_k=1)=p_k$.

\sk\noi Thinning the two-dimensional grid $\setZ^2$ is conceptually the same. 
Here the probability of keeping the point $(i,j)\in\setZ^2$ is $p_{i,j}$ where it is assumed that $\sum_{(i,j)\in \setZ^2}p_{i,j}<\infty$ whereby guaranteeing that the number of points of the thinned grid is finite. It suffices to apply bijectively $\setZ^2$ on $\setN$ by enumerating the points of $\setZ^2$ as $\{(i_n, j_n)\}_{n\geq 0}$, defining this bijection by $(i_n, j_n)\to n$. The rest is then obvious. 

\sk\noi This method may be useful when a sample of the germ point process is given (experimentally). The above thinning procedure can be adapted to this case. 

\sk\noi There is still an issue left aside in the presentation of the thinning procedure of the grid $\setN$. Can we really sample $T$? In fact one needs to have at disposition a closed expression of the distribution of this variable, in particular of the infinite product $\prod_{k\geq n+1}(1-p_k) $. If this is not possible, we may be lucky enough to find a dominating distribution function $q_n\geq p_n$ such $\sum_n q_n <\infty$ and such that the infinite product $\prod_{k\geq n+1}(1-q_k) $ is computable. One would then sample the thinned grid with thinning probability function $q_n$. A point of this dominating grid located at $k$ will the be kept with probability $p_k/q_k$ as a point of the desired sample.

\sk\noi For instance, try $q_n=1-e^{-\alpha_n}$ with $\sum_{n\geq 0} \alpha_n <\infty$ so that 
$$
\sum_{n\geq 0} q_n =\sum_{n\geq 0} 1-e^{-\alpha_n}  \leq \sum_{n\geq 0}\alpha_n <\infty \, . 
$$
The infinite products $\prod_{k\geq n+1}(1-q_k) $ should be computable, or equivalently, the sum $\sum_{n\geq 0} \alpha_n$ should be computable (and finite). This is the case for instance if $\alpha_n =C\frac{1}{n^2}$.

\sk\noi Note that the issue of computing the distribution of the number of points of the thinned point process is present even in the Poissonian case, where one needs to compute the integral $\int_{\setR^m} p(t) \, dt$. .

\sk\noi We now proceed to extension (b). 

\sk\noi It is a well-known fact that non-homogeneous Poisson process with intensity 
function $\lambda (t)$ can be obtained by projecting onto the time axis the points of an homogeneous Poisson 
process on $\setR^2$ of intensity $1$ 
which lie between the curve $y=\lambda (t)$ and the time axis.  
This is generalizable to point process admitting a stochastic intensity. 

\sk\noi The following result ([Grigelionis]) contains implicitly a simulation method for point processes with a stochastic intensity ([Ogata]). 

\sk\noi Let $(K, \K)$ be some measurable space. Given a history $\{{\cal F}_t\}_{t\in\setR}$, the point process $\overline{N}$ on $\setR\times K$ is called an ${\cal F}_t$-Poisson process if the following conditions are satisfied:

\sk (i) $\{{\cal F}_t\}_{t\in\setR}$ is a history of $\overline{N}$; 

\sk (ii) $\overline{N}$ is a Poisson process; and 

\sk (iii) for any $t\geq 0$, $S_t\overline{N}_+$ and ${\cal F}_t$ are independent ($S_t\overline{N}_+$ is the restriction of $N$ to $(t,\infty)$).

\sk\noi Let $\overline{N}$ be a $\F_t$-Poisson process on $\setR  \times \setR_+$ with intensity measure $dt  \times ds$. Let $f: \Omega \times \setR  \to \setR$ be a non-negative function that is $\P(\F_\cdot )$-measurable and such that the process 
$$
\lambda (t)\defeq f(t)
$$ 
is locally integrable. The point process $N$ defined by by 
$$
N(dt)\defeq \overline{N}((dt  \times [0, f(t)]) 
$$
admits the $\F_t$-stochastic intensityl $\lambda (t)$.

\sk\noi 
\subsection{Thinning a renewal point process} Let $N_0$ be an undelayed renewal sequence on $(0,+\infty)$: for $n\geq 1$, $X_{0,n}=S_1+\cdots +S_n$ where the sequence of non-negative random variables $\{S_n\}_{n\geq 1}$ is {\sc iid}, with a common distribution admitting a density $f$, with a corresponding failure rate $r(t)\defeq \frac{f(t)}{1-\int_0 f(s)\, ds}$ uniformly bounded by 
$M<\infty$. 

\sk\noi The stochastic intensity of such point process is $\lambda (t)=r(t-\theta_t)$ where $\theta_t$ is the position of the last point of $N_0$ that is $<t$. If we have at disposition a homogeneous Poisson process $\overline{N}$ on the strip $(0, +\infty)\times [0,M]$ with intensity $1$, the standard recursive procedure to generate $N_0$, based on the representation ($\star$), is the following. Given $X_{0,n}$, 
$X_{0,n+1}$ is smallest $t>X_n$ such that $\overline{N}$ has a point below the curve 
$y=r(t-X_{0,n}$. It then remains to thin this process with the thinning 
probability function $p(t)$ such that $\int_0^\infty p(t)\, dt<\infty$. 

\sk\noi To do this we shall do the thinning \emph{before} the construction of the basic renewal process. This is how. First construct a finite (due to the integrability condition on the thinning probability function) Poisson point process on the positive line of intensity $Mp(t)$, whose points are $t_1, \ldots, t_k$ in this order. Add to these points those of a Poisson process of intensity $M(1-p(t))$ to obtain a sequence $t'_1, t'_2, \ldots$ (there is an infinity of them but only those up to $t_k$ included will be used). Merge this sequence with the sequence $t_1, \ldots, t_k$. This merging produces a sample of a Poisson process of intensity $M$ on the time axis: $T_1, T_2, \ldots$. If $T_n\in \{t_1, \ldots, t_k\}$, set $X_n=1$, otherwise set $X_n=0$. Given $T_1, T_2, \ldots$, the sequence $\{X_n\}_{n\geq 1}$ is independent and the probability that $X_n=1$ is $p(T_n)$. 

\sk\noi Now let $\{W_n\}_{n\geq 1}$ be an {\sc iid} sequence uniformly distributed on $[0,1)$. The sequence $\{T_n,W_n\}_{n\geq 1}$ of points of $\setR^2$ form a Poisson process $\overline{N}$ of intensity $M$ on the strip $\setR_+\times (0,M)$. 
Construct the renewal germ point process: 
$$
N_0(dt)=\overline{N}(dt \times (0, r(t-\theta_t))
$$
where $\theta_t$ is the last point of $N_0$ that is $<t$ or $0$ if $N_0((0,t))=0$. 
The points of $N_0$ belong to the sequence $T_1, T_2, \ldots$. If a point of $N_0$ 
is $T_n$ keep it if and only if the corresponding $X_n$ is $1$. The surviving points are the points of a renewal process with failure rate $r(t)$ thinned with the probability function $p(t)$.

\sk\noi The case of a germ point process that is a delayed renewal point process is similar, {\it mutatis mutandis}. More generally, the case where the germ process is a point process with stochastic intensity with respect to its internal history $\lambda(t)\leq M$, receives a similar treatment, using the regenerative form of the stochastic intensity, as long as one is able to construct with finite computations a sample of the germ process. Here is an example where this is theoretically feasible.

\sk\noi 
\subsection{Exact sampling of a particular non-linear Hawkes point process} 
The so-called non-linear 
Hawkes process\index{Hawkes process!non-linear ---} (in its simplest form) is a 
point process $N$ on $\setR$ with the ${\cal F}^N_t$--intensity
$$
\lambda(t)\defeq \varphi\left(\int_{(-\infty,t)}h(t-s)N(ds)\right)
$$
where $\varphi:\setR\to\setR$ is a non-negative measurable function, and 
$h:\setR\to\setR$ is a measurable function (not necessarily non-negative) 
such that
$$
\pth{t<0\Rightarrow h(t)=0} \mbox{  and  } \int_{\setR_+}|h(t)|dt<\infty \, .
$$

\sk\noi Suppose that $h:\setR\to\setR$ has a bounded support $[0,a]$ in the sense that 
$h(t,z)\neq 0$ implies that $t\in [0,a]$. Suppose in addition that 
$\varphi$  is bounded (say, 
by $\Lambda<\infty$). 

\sk\noi A construction of a stationary version of this point process is as follows. Let $\tilde N^a$ be the 
point 
process formed by the points $\tilde T_n$ of a Poisson process $\tilde N$ of intensity $\Lambda$ such that 
$\tilde T_n-\tilde T_{n-1}>a$, and call $\{\tilde T^a_n\}_{n\in\setZ}$ the 
sequence of points of $\tilde N^a$. These are ``regeneration points'', because for all $t\in[\tilde T^a_n,\tilde T^a_{n+1})$
$$
\lambda(t)=\varphi\left(\int_{(\tilde T^a_n,t)}h(t-s)N(ds)\right)
$$
does not depend on $N$ before $\tilde T^a_n$. Thus we have an explicit form for $\lambda(t)$ for all $t$ that 
does not require knowledge of the whole past of $N$.

\sk\noi We shall not give a theoretical description of extension (c), since the following example clearly gives the method, which applies each time that the point process to be thinned is constructed ``below a Poisson process''.

\subsection{Thinning the Mat\'ern hard-core model} 
([Mat\'ern]) Let $N$ be a homogeneous Poisson process on $\setR^m$ with intensity $\lambda$. Let $\{X_n\}_{n\in\setN}$ be its sequence of points. The Mat\'ern model is a point process $\widetilde{N}$ obtained by thinning $N$ in such a way that all pairs of points of the thinned point process are at least at a distance $r>0$ apart. The thinning is done as follows. Let $\{U_n\}_{n\in\setN}$ be an {\sc iid} sequence of real random variables uniformly distributed on the interval $[0,1]$, independent of $N$. A point $X_n$ of $N$ is retained as a point of $\widetilde{N}$ if and only if 
$$ 
U_n < U_k \mbox{ for all } k\neq n \mbox{ such that } X_k\in \overline{B}(X_n; r) \, ,
$$
where $\overline{B}(x;r)$ denotes the closed ball of center $x$ and radius $r$. 

\sk\noi The ``thinning first'' method applies. First generate the (finite) thinned point process, to obtain a Poisson process $N_1$ of intensity $\lambda p(x)$. Then generate the points of a Poisson process $N_2$ with intensity $\lambda (1-p(x))$. Only a finite number of those points have to be generated, those who lie at distance less that $r$ from the thinned Poisson process. Superposing them to those of the thinned Poisson process, one obtains a Poisson point process$N=N_1+N_2$, 
with intensity $\lambda$, in a limited region. The points of $N_1$ are then thinned according to Mat\`ern's construction, using the points of $N$. The surviving points of $N_1$ form a Mat\`ern point process thinned according to the retaining probability function $p(x)$.

\section{Exact sampling of Hawkes processes}

One seeks to obtain on $W\defeq [0,a]$ an exact sample of a linear Hawkes process $N$ on the line with random fertility rate $h(t,Z)$, where $Z$ is a random element in some measurable space $(K, \K)$, and such that $\rho\defeq E\cro{\int_0^\infty h(t,Z)\, dt} <1$. This process is a cluster point process  
where the germ point process $N_0$ is a Poisson process with intensity 
function $\widetilde{\mu}$ and where for each $n\geq 1$, $Z_n$ is a branching point process of random fertility rate $h(t,Z)$  with a single ancestor located at $0$. The following result will be needed. 
Consider the univariate branching Hawkes process on the line with a single ancestor at the origin of times, and with random fertility rate $h(t,Z)$ such that 
$$
t<0 \rightarrow h(t,z)=0 \mbox{ for all } z\in K 
$$
and 
$$
\int_{0}^{\infty}E[h(t,Z_1)]dt<1 \, .
$$
The shifted process $S_{t}N$ 
converges in distribution to the empty process, and if moreover
$$
\int_{0}^{\infty}t E[h(t,Z_1)]dt<\infty
$$
the convergence is in variation (Br\'emaud and Massouli\'e, 1996). Therefore there exists a finite non-negative random variable $L$, called the extinction time, such that $N$ is empty on $(L,+\infty)$.

\sk\noi In principle, the Brix--Kendall perfect simulation method of $N$ on $W=[0,a]$ applies. It consists in two steps. First, one generates a Poisson point process on $(-\infty, 0]$ with intensity $\widetilde{\mu}(t) P(L\geq -t)$ where $L$ is the typical length (extinction time) of a branching point process on $\setR_+$ with random fertility rate $h(t,Z)$ and with a single ancestor located at $0$, and then, for each of the points of this Poisson process, one generates samples of a branching point process of random fertility rate $h(t,Z)$  with a single ancestor located at $\widetilde{T}_n$ until one obtains a sample with at least one point in $(0,\infty)$.  
(The rationale is that the point of $N_0$ located at $-t$ has a probability $P(L\geq -t)$ of having points in $[0,\infty)$\,\footnote{Of course $W=[0,a]\subseteq [0,\infty)$, but a little thought will convince the reader that we can pretend that the goal is to sample $N$ on $[0,\infty)$ without additional cost since the points of $N_0$ after $a$ will not be used.}. 

\sk\noi Exact sampling requires that the Poisson process on $(-\infty, 0]$ of intensity $\widetilde{\mu}(t) P(L\geq -t)$ finite. This implies some restrictions. For instance, if $\widetilde{\mu}$ is bounded, the condition $\int_0^\infty P(L>t)\, dt =E[L]<\infty$ will guarantee that. In particular, $L$ must be finite, that is convergence of the cluster to the empty point process (the stationary state) must take place in finite time. The convergence is then in variation, and a sufficient condition for this is $E\cro{\int_0^\infty t h(t,Z)\, dt}<\infty$.

\sk\noi Application of the Brix--Kendall method requires the construction of a Poisson process on $\setR_+$\,\footnote{Here, to facilitate notation, time is reversed, so that the surviving points appear to be on the positive line rather than on the negative line.} of intensity $\widetilde{\mu}(t) P(L > t)$ but the exact form of $P(L > t)$ is not known. However, suppose that we know explicitly sequences of non-negative functions $\{l_n\}_{n\geq 1}$ and $\{u_n\}_{n\geq 1}$ respectively non-decreasing and non-increasing, and both converging pointwise to $P(L > t)$ in such a way that $||u_n-l_n||_\infty \defeq \sup_{t\geq 0} |u_n(t)-l_n(t)|\, dt \to 0$. Suppose moreover that $\int_0^t u_0(t)\widetilde{\mu}(t) \, dt <\infty$. This is the case in particular if $u_0(t)=1-G(t)$ for some cumulative distribution function on $\setR_+$ of finite mean and $\sup_{t\geq 0}\widetilde{\mu}(t)<\infty$. 

\sk\noi The construction goes as follows. First generate the points (in finite number) of a Poisson process of intensity $u_0(t)\widetilde{\mu}(t)$. Let $t_1$, \ldots, $t_k$ be these points. Generate an {\sc iid} sequence $V_1$, \ldots, $V_k$, of random variables uniformly distributed on $[0,1]$. Under the curve $y=u_0(t)\widetilde{\mu}(t)$, the points $(t_1, V_1u_0(t_1)\widetilde{\mu}(t_1))$, \ldots, $(t_k, V_1u_0(t_k)\widetilde{\mu}(t_k))$ of $\setR_+^2$ form a Poisson process of intensity $1$. Note that the probability that any of these points lie on the curve $y=\widetilde{\mu}(t) P(L > t)$ is null, and therefore, it will happen for a finite $n$ that none of the points $(t_1, V_1u_0(t_1)\widetilde{\mu}(t_1))$, \ldots, $(t_k, V_1u_0(t_k)\widetilde{\mu}(t_k))$ lies between the curves 
$y=u_n(t)\widetilde{\mu}(t)$ and $y=l_n(t)\widetilde{\mu}(t)$. For the simulation, keep only the points $t_i$ such that $V_iu_0(t_i)\widetilde{\mu}(t_i)<\ell_n(t_i)$ since these points are exactly those lying under the curve $y=
P(L>t)$.

\sk\noi It remains to find the approximating functions $u_n$ and $l_n$. 
Let $\hat{N}_k$ be $\hat{N}$ restricted to the generations $0$, $1$, \ldots, $k$. In particular, 
$\lim_{k\uparrow \infty} \hat{N}_k((t,\infty))= \hat{N}((t,\infty))$. Compute 
$E\cro{e^{-\theta \hat{N}_k((t,\infty))}}$ for any $\theta >0$, take the limit as $k\uparrow \infty$ to obtain $E\cro{e^{-\theta \hat{N}((t,\infty))}}$, and then obtain the result from the remark that  
$$
\lim_{\theta\uparrow +\infty}E\cro{e^{-\theta \hat{N}((t,\infty))}}=P(\hat{N}((t,\infty))=0)\, .
$$  
To compute $E\cro{e^{-\theta \hat{N}_k((t,\infty))}}$, introduce the random elements $Y_k$ with values in $(M_p(\setR_+), \M_p(\setR_+))$ whose distribution is that of $N_k$. Let now $\{Y_{k, n}\}_{n\in\setN}$ be an {\sc iid} sequence of random elements distributed as $Y_k$. The distribution of $\hat{N}_k((t,\infty))$ is the same as that of 
$$
\sum_{n\geq 0}1_{\{\tau_n\leq t\}} Y_{k-1,n} (t-\tau_n)
$$
where the $\tau_n$'s are the points of the first generation of $\hat{N}$, forming a Cox point process of random intensity $h(t,Z)$. Therefore, 
\begin{align*}
E\cro{e^{-\theta \hat{N}_k((t,\infty))}}
&=E\cro{e^{-\theta \sum_{n\geq 0}1_{\{\tau_n\leq t\}} Y_{k-1,n} (t-\tau_n)}} \\
&=E\cro{e^{-\theta \sum_{n\geq 0}1_{\{\tau_n\leq t\}} f(Y_{k-1,n},t-\tau_n)}}
\end{align*}
where for $\mu\in M_p(\setR_+)$ and $u\in\setR_+$, $f(\mu, u)\defeq \mu((u,\infty))$, so that 
\begin{align*}
&E\cro{e^{-\theta \sum_{n\geq 0}1_{\{\tau_n\leq t\}} f(Y_{k-1,n},t-\tau_n)}}\\
&\hphantom{aaa}E\cro{\exp\acc{\int_0^t \int_{M_p(\setR_+)}
\pth{e^{-\theta f(\mu, t-s)}-1} h(s,Z)\, Q_{k-1}(d\mu)   
}}\\
&\hphantom{aaa}E\cro{\exp\acc{\int_0^t \int_{M_p(\setR_+)}\pth{e^{-\theta \mu((t-s,\infty))}-1} h(s,Z)\, Q_{k-1}(d\mu) ds  }}\, .
\end{align*} 
Now, 
$$
\int_{M_p(\setR_+)}\pth{e^{-\theta \mu((t-s,\infty))}-1}Q_{k-1}(d\mu)=E\cro{ e^{-\theta \hat{N}_{k-1}((t-s,\infty))}}
$$
and therefore 
$$
E\cro{e^{-\theta \hat{N}_k((t,\infty))}}=E\cro{\exp\acc{\int_0^t \pth{E\cro{ e^{-\theta \hat{N}_{k-1}((t-s,\infty))}}-1      }h(s,Z)\, ds} }\, .
$$
Taking the limit as $k\uparrow \infty$ yields 
$$
E\cro{e^{-\theta \hat{N}((t,\infty))}}=E\cro{\exp\acc{\int_0^t \pth{E\cro{ e^{-\theta \hat{N}((t-s,\infty))}}-1      }h(s,Z)\, ds} }
$$
and taking the limit as $\theta\to +\infty$ gives, with $f(t)=P(\hat{N}((t,\infty))>0  )$,   
$$
f(t)=E\cro{\exp\acc{-\nu(t,Z)+\int_0^t f(t-s) h(s,Z) \, ds}} \, ,
$$
where $\nu(t,Z)\defeq \int_0^t h(s,Z)\, ds$. The right-hand side will be denoted by $\Phi(f)$. We are therefore concerned with the equation 
$$
f=\Phi(f)\, , \, (f\in \A),   
$$
where $\A\defeq \{f: (\setR_+, \B(\setR_+)) \to ([0,1]), \B([0,1]))\}$. We show that the solution $F(t)\defeq P(\hat{N}((t,\infty))>0$ is the unique solution, and that there exists sequences 
$\{g_n\}_{n\geq 1}$ and $\{h_n\}_{n\geq 1}$, respectively non-decreasing and non-increasing, with the common limit $F$ and such that    
$$
||f_n-g_n||_\infty \to 0\, . 
$$

\sk\noi We summarize the main steps of the approximation in [M{\o}ller and Rasmussen]

\sk\noi (i) The sequence of functions $\{f_n\}_{n\geq 1}$ defined by $f_n(t)=P(\hat{N}_n((t,\infty))>0  )$ is non-decreasing and has $F$ for limit and that it satisfies the recurrence 
$$
f_{n}=\Phi (f_{n-1})\, , \, n\geq 1 \, ,
$$ 
with $f_0(t)\equiv 1$.

\sk\noi (ii) Let $f\in\A$ and let $\Phi^{(n)}$ be defined recursively by $\Phi^{(0)}(f)\defeq f$, 
$\Phi^{(n)}(f)\defeq\Phi^{(n-1)}(\Phi(f))$, and let $f_n\defeq\Phi^{(n)}(f)$. 
We have that 
\begin{align*} 
f&\leq g \Rightarrow f_n\leq g_n \, , \\
f&\leq \Phi(f)\Rightarrow \{f_n\}_{n\geq 1} \mbox{ is non-decreasing }\, ,\\
f&\geq \Phi(f)\Rightarrow \{f_n\}_{n\geq 1} \mbox{ is non-increasing } \, .
\end{align*}

\sk\noi (iii) $\Phi$ is a contraction on $\A$ with respect to the sup norm. More precisely  
$$
f,g\in\A \Rightarrow ||\Phi(f)-\Phi(g) ||_\infty \leq \rho ||f-g ||_\infty \, .
$$

\sk\noi (iv) $F$ is the unique fixed point of $\Phi$.

\sk\noi (vi) $||f_n-F ||_\infty\leq \frac{\rho^n}{1-\rho}||\Phi(f)-f ||_\infty$.

\sk\noi (vii) If $f\leq \Phi(f)$ or $f\geq \Phi(f)$, then $f_n\to F$ from below or from above respectively. This follows from (ii).

\sk\noi The functions $u_n$ and $\ell_n$ are now defined. Since $0\leq \Phi(0)$, $f_n^u \defeq \Phi^{(n)}  (0)\downarrow F$. Take $\ell_n(t)=1-f_n^u(t)$ so that $\ell_n\uparrow 1-F$. M{\o}ller and Rasmussen exhibit a cumulative distribution function $G$ with finite mean such that $G\geq \Phi(G)$. Then $f_n^\ell \defeq \Phi^{(n)}  (G)\uparrow F$. Take $u_n(t)=1-f_n^\ell(t)$ so that $u_n\downarrow 1-F$. This is crucial since the algorithm starts by constructing (sampling) a Poisson process on $\setR_+$ of intensity $u_0(t)\widetilde{\mu}(t)=(1-G(t))\widetilde{\mu}(t)$ of finite mass\,\footnote{The choice $f_n^\ell \defeq \Phi^{(n)}  (1)$ would lead to the correct limit $f_n^\ell\uparrow F$, but in this case the Poisson process of intensity 
$u_0(t)\widetilde{\mu}(t)=\widetilde{\mu}(t)$ could have an infinite number of points, for instance if $\widetilde{\mu}(t)$ is a positive constant.}.

\sk\noi The M{\o}ller--Rasmussen thinning procedure is adaptable to the case of a germ process that is a renewal process, or a delayed renewal point process, and more generally, a point process with stochastic intensity with respect to its internal history $\lambda(t)\leq M$. The case of a point process on the linear grid, or of a Mat\`ern point process on the line, are also amenable to exact sampling once the thinning of the germ point process is feasible. It suffices to proceed as in the standard case because the M{\o}ller--Rasmussen algorithm is able to construct the thinned Poisson process as well as the Poisson process over he line $y=Mp(t)$. The details are left for the reader. 

\section{Conclusion and summary} 

This article gives a formal proof of the validity of the Brix--Kendall exact sampling algorithm via Laplace transforms. This proof is the basis for several extensions, in particular to a germ process that is a grid, or a subset of a given deterministic point configuration (for example, in a communications context, antenna locations), or, on the line, a renewal process with a bounded fertility rate. The case of unbounded fertility rate will follow the same basic idea and presents no additional conceptual difficulty. In general, any germ point process that can be constructed ``under a Poisson process'', such as the Mat\`ern hard-core point process, fits the proposed extension. The article also shows that the Brix--Kendall exact sampling method applies in a natural way to the exact sampling of Boolean models. The M{\o}ller--Rasmussen exact sampling algorithm of Hawkes point processes on the line with light-tail fertility rate has been extended to the situation where the germ point processes is a renewal process with bounded failure rate, or some point process with a stochastic intensity.

\section*{References}

\sk\noi P. Br\'emaud, {\it Point Processes}, to appear, 2017.

\sk\noi P. Br\'emaud and L. Massouli\'e, ``Stability of non--linear 
Hawkes processes'', {\it Annals of Probability\/} {\bf 24} 
(3) 1563-1588 (1996).

\sk\noi A. Brix and W.S. Kendall, ''Simulation of cluster point processes without edge effects'', {\it Adv. Appl. Prob.}, {\bf 34}, 267-280 (2002).

\sk\noi D.J. Daley, D. Vere--Jones, {\it An Introduction to the 
Theory of Point Processes\/}, Springer, NY (1988, 2003).

\sk\noi B. Grigelionis, ``On the representation of integer-valued random measures by means of stochastic integrals with respect to a Poisson measure'', {\it Litovsk. Mat. Sb.} {\bf 11}, 93-108 (in Russian) (1971).

\sk\noi A.G. Hawkes and D. Oakes ``A cluster process representation 
of a self--exiciting point process'', {\it J. Appl. Proba.\/} {\bf 11}, 
493--503 (1974).

\sk\noi B. Mat\'ern, ``Spatial Variation'', {\it Meddelandenfran Statens Skogsforskningsinstitut}, {\bf 49} (5), 1-44 (1960).

\sk\noi J. M{\o}ller and J.G. Rasmussen. ``Perfect simulation of Hawkes processes'', Advances in Applied Probability, 37, 629--646 (2004).

\sk\noi Y. Ogata, ``On Lewis’ simulation method for point processes'', 
{\it IEEE Trans. Inf. Theory}, {\bf 27}, 23–31 (1981).

\end{document}